\documentclass[11pt,a4paper]{article}
\usepackage{amsfonts}
\usepackage{}

\usepackage{graphicx}
\usepackage{amsfonts,amsmath, amssymb}

\newtheorem{theo}{\textbf{\ \ \quad Theorem}}[section]

\newtheorem{lem}{\textbf{\ \ \quad Lemma}}[section]
\newtheorem{remark}{\textbf{\ \ \quad Remark}}[section]

\newtheorem{prop}{\textbf{\ \ \quad Proposition}}[section]
\newtheorem{defi}{\textbf{\ \ \quad Definition}}[section]

\newcommand{\lbl}[1]{\label{#1}}

\newcommand{\be}{\begin{equation}}
\newcommand{\ee}{\end{equation}}
\newcommand\bes{\begin{eqnarray}}
\newcommand\ees{\end{eqnarray}}
\newcommand{\bess}{\begin{eqnarray*}}
\newcommand{\eess}{\end{eqnarray*}}

\newcommand{\nm}{\nonumber}

\newcommand{\ds}{\displaystyle}

{\setlength\arraycolsep{2pt}}

\title{Noise and Stability in Reaction-diffusion Equations}

\author{Guangying Lv$^{a}$,  Jinlong Wei$^b$, Guang-an Zou$^{c}$\\
\\   {\small \it $^a$College of Mathematics and Statistics,
Nanjing University of Information}\\
{\small \it  Science and Technology, Nanjing 210044, China}\\
{\small \tt gylvmaths@126.com}\\
{\small \it $^b$School of Statistics and Mathematics, Zhongnan University of}\\ {\small \it Economics and Law, Wuhan 430073, China}\\ {\small \tt  weijinlong.hust@gmail.com }\\
{\small \it $^c$Institute of Applied Mathematics, Henan University}\\ {\small \it Kaifeng 475001, China}
\\ {\small \tt zouguangan@henu.edu.cn}}

\begin{document}
\maketitle

\begin{abstract}
We study the stability of reaction-diffusion equations in presence of noise. The relationship of stability of solutions between the stochastic ordinary different equations and the
corresponding stochastic reaction-diffusion equation is firstly established. Then, by using the Lyapunov method, sufficient conditions for mean square and stochastic
stability are given. The results show that the multiplicative noise can make the solution stable, but the
additive noise will be not.



\end{abstract}

\baselineskip=15pt

\section{Introduction}
\setcounter{equation}{0}

The stability of solutions  is an important issue in
the theory of PDEs (partial differential equations), which has been studied by
many authors \cite{YLbook}. There are a lot of sufficient conditions to assure that the
solutions are stable or unstable. We note that
noise always exists in the real world.
The reasons may be  that the parameter is obtained by different measurement,
and we can not get the real value, so we consider the ordinary (partial) differential
equations with noise perturbation is available.
In other words, in microscopic world, to describe the particle moving law must
be stochastic ordinary (partial) differential equations. In macroscopic world,
we often consider the case that the coefficient in equations is random or stochastic.
 Usually, if we consider the
role of noise, we have two cases. First case: the noise is regarded as a small perturbation.
In this case, the structure of solutions will not be changed and the biggest possible change is
long-time behavior of the solutions, that is to say, the conditions of stability or un-stability
may be different from the deterministic case. Second case: the noise is strong, such as $u^\gamma dW_t$,
where $\gamma>1$, $u$ is the unknown function and $W_t$ is the noise. In this case, the
structure of solutions will  be changed. More precisely, the noise can induce the solutions
blow up in finite time.

When the noise appears in a deterministic PDE, the impact of noise on solutions will
be the first thing to be considered. In the present paper, we aim to study impact of noise
on stability of solutions.

Initially, we recall some known results about the impact of noise.
Flandoli et al. \cite{FF2012,FGP2010} proved that the noise can make
the transport equations well-posedness and the noise can prevent the
singularities in linear transport equations. Chow \cite{C2011} obtained that the noise can induce
singularities (finite time blow up of solutions), also see \cite{LD2015,LW2020}.
There are a lot of work about the impact of noise on different PDEs, for example,
Hamilton-Jacobi equations \cite{GG2019},
conservation law \cite{GS2017,GL2019}, porous media equations \cite{DGG2019},
quasilinear degenerate parabolic-hyperbolic equations \cite{GH2018}, and so on.
Besides, the impact of noise on regularity of solutions of parabolic equations
has been studied by Lv et al. \cite{LGWW2019}.

Forty years ago, there are a lot of works about the impact of noise. The main issue
is that the noise can stabilize the solution of ordinary differential equations,
see the book \cite{Kh2011,Maobook1991,Maobook1994}. Meanwhile, the stochastic stability of
functional differential equations is also considered by Mackey-Nechaeva \cite{MN1994}.
In the book \cite{DZbook}, the long time behavior of solutions was considered in Chapter 11 and
sufficient condition of mean square stable is given, see Theorem 11.14. More precisely, Da Prato- Zabczyk
studied the following equation
     \bes\left\{\begin{array}{lll}
   dX=AXdt+B(X)dW_t, \\[1.5mm]
   X(0)=x\in H,
    \end{array}\right.\lbl{1.1}\ees
where $H$ is a Hilbert space, $A$ generates a $C_0$ semigroup and $B\in L(H;L_2^0)$ (see
p309 of \cite{DZbook} for more details). They proved that the following statements are equivalent

(i) There exists $M>0$, $\gamma>0$ such that
   \bess
\mathbb{E}|X(t,x)|^2\leq Me^{-\gamma t}|x|^2, \ \  \ t\geq0.
    \eess

(ii) For any $x\in H$ we have
   \bess
\mathbb{E}\int_0^\infty|X(t,x)|^2dt<\infty.
   \eess
It is remarked that the nonlinear term is not considered in the book \cite{DZbook}. Moreover,
the example in \cite{DZbook} is the stochastic reaction-diffusion equation on the bounded domain.
Liu-Mao \cite{LM1998} also considered the stability of trivial solution $0$ on the bounded domain.
Wang-Li \cite{WL2019} considered the
stability and moment boundedness of the stochastic linear age-structured model.
Although in the book \cite{Cb2007}, Chow gave a abstract result to study the stability of null solution
(see page 233), the concrete form was not given.
In the present paper, we consider the concrete model and generalized the classical results.
What's more, we find some difference between stochastic partial differential equations (SPDEs)
and stochastic differential equations (SDEs).

We discuss the impact of different kinds of noise on stability.
Some interesting results are obtained: the additive noise will have a "bad" effect
and some multiplicative noise has "good" effect. In other words, the multiplicative noise can make the solution stable, but the
additive noise will not, which is new for SPDEs.
This is different from SDEs, see Remark \ref{r3.2}.
What's more, we obtain a new result for stability theory of SDEs, see Theorem \ref{t3.6}.

This paper is arranged as follows.   In Sections 2,
we will present some preliminaries. Sections 3 and 4 are concerned with the bounded domain and
the whole space, respectively.

\section{Preliminaries}
\setcounter{equation}{0}
In the present paper, we always assume $W_t$ is a one-dimensional standard Wiener process and $W_t(x)$ is a Wiener random field which are defined on a complete probability
space $(\Omega,\mathcal {F},\mathbb{P})$. Firstly, we recall the definitions of stochastic stability.
We only consider the stability of constant equilibrium of stochastic
reaction-diffusion equations.
Consider the following stochastic reaction-diffusion equations
  \bes\left\{\begin{array}{lll}
   du=(\Delta u+f(u))dt+\sigma(u)dW, \ \ \qquad t>0,&x\in  \mathbb{R}^d,\\[1.5mm]
   u(x,0)=u_0(x), \ \ \ &x\in \mathbb{R}^d,
    \end{array}\right.\lbl{2.2}\ees
where $W=W_t$ or $W_t(x)$. The Wiener random field
can be chosen to have the following properties:
$\mathbb{E}W_t(x)=0$ and its covariance function $q(x,y)$ is given by
    \bess
\mathbb{E}W_t(x)W_s(y)=(t\wedge s)q(x,y), \ \ \ x,y\in\mathbb{R}^d,
   \eess
where $t\wedge s=\min\{t,s\}$ for $0\leq t,s\leq T$, or can be chosen as one-dimension Brownian motion.
Without loss of generality, we suppose that $f(0)=\sigma(0)=0$, that is to say, $0$ is a trivial
solution to the first equation of problem (\ref{2.2}). Let $\|\cdot\|$ denote
as some norm with respect to the spatial variable.

\begin{defi} \lbl{d2.1} The trivial solution $0$ is called mean square stable if for any
$\varepsilon>0$, there exists $\delta=\delta(\varepsilon)>0$ such that for any initial
data $u_0$,
   \bess
\|u_0\|\leq\delta\ \ \ {\rm implies}\ \ \ \mathbb{E}\|u(\cdot,t)\|^2<\varepsilon, \ \forall \ t\geq0,
     \eess
and exponentially mean square stable, if there exists two positive constants $c_1$ and $c_2$ such that
     \bess
\mathbb{E}\|u(\cdot,t)\|^2<c_1e^{-c_2t}\|u_0\|^2, \ \forall \ t\geq0.
     \eess
\end{defi}

\begin{defi} \lbl{d2.2} The trivial solution $0$ is called stochastically stable if for any
$\varepsilon_1>0$ and $\varepsilon_2>0$, there exists $\delta=\delta(\varepsilon_1,\varepsilon_2)>0$ such that for $t>0$
the solution $u$ satisfies
   \bess
\mathbb{P}\left\{\sup_{t>0}\|u(\cdot,t)\|\leq\varepsilon_1\right\}\geq1-\varepsilon_2 \ \ \ {\rm for}\ \ \ \|u_0\|\leq\delta.
     \eess
\end{defi}

Before ending this section, we consider the relationship of the stability of solutions between
stochastic differential equations and stochastic reaction-diffusion equations.

Assume that $u=0$ is a trivial solution to the first equation of problem (\ref{2.2}), then
$u=0$ will be a trivial solution to the following equation
   \bes
du=f(u)dt+\sigma(u)dW_t.
   \lbl{2.3}\ees
Problem (\ref{2.3}) means
equation (\ref{2.3}) with the initial data $u_0$ (independent of $x$). Initially, we introduce a definition for (\ref{2.3}).

\begin{defi} \lbl{d2.1} The trivial solution $0$ is called mean square stable for problem (\ref{2.3}) if for any
$\varepsilon>0$, there exists $\delta=\delta(\varepsilon)>0$ such that for any initial
data $u_0$,
   \bess
|u_0|\leq\delta\ \ \ {\rm implies}\ \ \ \mathbb{E}|u(\cdot,t)|^2<\varepsilon, \ \forall \ t\geq0,
     \eess
and exponentially mean square stable, if there exists two positive constants $c_1$ and $c_2$ such that
     \bess
\mathbb{E}|u(\cdot,t)|^2<c_1e^{-c_2t}|u_0|^2, \ \forall \ t\geq0,
     \eess
and stochastically stable if for any
$\varepsilon_1>0$ and $\varepsilon_2>0$, there exists $\delta=\delta(\varepsilon_1,\varepsilon_2)>0$ such that for $t>0$
the solution $u$ satisfies
   \bess
\mathbb{P}\left\{\sup_{t>0}|u(\cdot,t)|\leq\varepsilon_1\right\}\geq1-\varepsilon_2 \ \ \ {\rm for}\ \ \ |u_0|\leq\delta.
     \eess
\end{defi}
In order to establish the
relationship between problems (\ref{2.2}) and (\ref{2.3}), we need the following lemma.
Let $\eta(r)=r^-$ denote the negative part of $r$ for $r\in\mathbb{R}$.
Set
   \bess
k(r)=\eta^2(r),
   \eess
so that $k(r)=0$ for $r\geq0$ and $k(r)=r^2$ for $r<0$. For
$\epsilon>0$, let $k_\epsilon(r)$ be a $C^2$-regularization of $k(r)$ defined by
   \bess
k_\epsilon(r)=\left\{\begin{array}{lllll}
r^2-\ds\frac{\epsilon^2}{6}, \ \qquad \ &r<-\epsilon,\\[1mm]
-\ds\frac{r^3}{\epsilon}\left(\ds\frac{r}{2\epsilon}+\frac{4}{3}\right),\ \ &-\epsilon\leq r<0,\\[1mm]
0, \ \ \ &r\geq0.
   \end{array}\right.\eess
Then one can check that $k_\epsilon(r)$ has the following properties.
  \begin{lem}\lbl{l2.1}{\rm\cite[Lemma 3.1]{LD2015}} The first two derivatives $k'_\epsilon,\,k''_\epsilon$
  of $k_\epsilon$ are continuous and satisfy the conditions: $k'_\epsilon(r)=0$ for $r\geq0$;
  $k'_\epsilon\leq0$ and $k''_\epsilon\geq0$ for any $r\in\mathbb{R}$. Moreover, as $\epsilon\rightarrow0$,
  we have
    \bess
k_\epsilon(r)\rightarrow k(r) \ \ {\rm and} \ \ k'_\epsilon(r)\rightarrow-2\eta(r),
   \eess
and the convergence is uniform for $r\in\mathbb{R}$.
   \end{lem}

Under the condition that
both problems (\ref{2.2}) and (\ref{2.3}) have a unique strong solution, we have the following result.
Here we focus on the stability of solutions and we do not talk about the
existence of solutions.
   \begin{theo}\lbl{t2.1} Assume that $W_t$ is a one-dimension Wiener process, $f$ and $\sigma$ satisfy the global Lipschitz condition,
and $f(0)=\sigma(0)=0$. Then $0$ is a stable (exponentially mean square stable or stochastically stable)
trivial solution  of (\ref{2.2}) (with the $L^\infty$ norm in spatial variable) if and only if
$0$ is a stable (exponentially mean square sable or stochastically stable) trivial solution of (\ref{2.3}).
   \end{theo}

{\bf Proof.} Clearly, if $0$ is a stable (exponentially mean square stable or stochastically stable) trivial solution of (\ref{2.2}),
then $0$ is also a stable trivial solution of (\ref{2.3}) since we can choose the initial data $u_0$ which is spatial variable independent.

Conversely, we will check if if $0$ is a stable (exponentially mean square stable or stochastically stable) trivial solution of (\ref{2.3}),
then $0$ is a stable trivial solution of (\ref{2.2}) as well.

Assume $0$ is a stable (exponentially mean square stable or stochastically stable) trivial solution of (\ref{2.3}),
then for any $\varepsilon>0$, there exists constant $\delta=\delta(\varepsilon)>0$ such that
   \bess
\mathbb{E}|u_\pm(t)|^2<\varepsilon,
   \eess
where $u_\pm(t)$ is the solution of (\ref{2.3}) with initial data $\pm2\delta$.

Let $u(x,t)$ be the unique solution of (\ref{2.2}) with initial data $u_0(x)$ which satisfies $|u_0(x)|\leq\delta$.  For any fixed $T>0$, we will prove that
   \bes
u_-(t)\leq u(x,t)\leq u_+(t), \ \ \forall \ (x,t)\in\mathbb{R}^d\times[0,T], \ a.s..
   \lbl{2.4}\ees
If the above inequality holds, then we get the desired result.
In order to the inequality (\ref{2.4}), we let
 $v=u(x,t)-u_-(t)$, then $v$ satisfies
  \bes\left\{\begin{array}{lll}
   dv=(\Delta v+f(u)-f(u_-))dt+(\sigma(u)-\sigma(u_-))dW_t, \ \ \qquad t>0,&x\in  \mathbb{R}^d,\\[1.5mm]
   v(x,0)=u_0(x)+\delta\geq0, \ \ \ &x\in \mathbb{R}^d.
    \end{array}\right.\lbl{2.5}\ees
Note that $\nabla u$ makes sense in $L^\infty(\mathbb{R}^d)$, thus we can not
take $\Phi_\varepsilon(v(t))=(1,k_\epsilon(v(t)))$, which is different from
those in \cite{C2011,LW2020}. We need introduce a new test function.

For any $R>0$, $B_R(0)$ denotes a ball centered in $0$ with radius $R$. Let $\phi_1$
be the eigenvalue function of Laplacian operator on $B_R(0)$ with respect to
the first eigenvalue $\lambda_1$, i.e.,
  \bess\left\{\begin{array}{llll}
-\Delta \phi_1=\lambda_1 \phi_1, \ \ \ \ \ \ \  \ {\rm in} \  B_R(0),\\
\ \phi_1=0, \ \ \qquad\ \ \ \qquad  {\rm on}\ \partial B_R(0).
   \end{array}\right. \eess
Denote $\psi\in C^2(\mathbb{R}^d)$ satisfying
    \bess
\psi(x)=\left\{\begin{array}{lll}
\phi_1(x), \ \ \ & x\in\bar B_R(0),\\
0,\ \ \ & x\in\mathbb{R}^d\setminus B_R(0).
   \end{array}\right.\eess
Define
   \bess
\Phi_\varepsilon(v(t))=(\psi,k_\epsilon(v(t)))=\int_{\mathbb{R}^d}\psi(x) k_\epsilon(v(x,t))dx.
   \eess
By It\^{o}'s formula, we have
    \bess
\Phi_\varepsilon(v(t))&=&\Phi_\varepsilon(v_0)+\int_0^t\int_{\mathbb{R}^d}\psi(x) k_\epsilon'(v(x,s))\Delta v(x,s)dxds\\
&&+\int_0^t\int_ {\mathbb{R}^d}\psi(x) k_\epsilon'(v(x,s))(f(u(x,s))-f(u_-(s)))dxds\\
&&+\int_0^t\int_ {\mathbb{R}^d}\psi(x) k_\epsilon'(v(x,s))(\sigma(u(x,s))-\sigma(u_-(s)))dxdW_s\\
&&+\frac{1}{2}\int_0^t\int_{\mathbb{R}^d}\psi(x)k_\epsilon''(v(x,s))(\sigma(u(x,s))-\sigma(u_-(s)))^2dxds.
    \eess
By using the facts $k_\epsilon''\geq0,\,\psi(x)\geq0$, we have
   \bess
&&\int_{\mathbb{R}^d}\psi(x) k_\epsilon'(v(x,s))\Delta v(x,s)dx\\
&=&-\int_{\mathbb{R}^d}\psi(x) k_\epsilon''(v(x,s))|\nabla v(x,s)|^2dx- \int_{\mathbb{R}^d}k_\epsilon'(v(x,s))\nabla\psi(x)\cdot\nabla v(x,s)dx\\
&\leq&- \int_{\mathbb{R}^d}\nabla k_\epsilon(v(x,s))\cdot\nabla\psi(x)dx\\
&=&\int_{\mathbb{R}^d} k_\epsilon(v(x,s))\Delta\psi(x)dx\\
&=&-\lambda_1\int_{\mathbb{R}^d} k_\epsilon(v(x,s))\psi(x)dx\leq0.
   \eess
Consequently,
     \bess
\Phi_\varepsilon(v(t))
&\leq&\Phi_\varepsilon(v_0)+\int_0^t\int_{\mathbb{R}^d}\psi(x)k_\epsilon''(v(x,s))\left(\frac{1}{2}(\sigma(u(x,s))-\sigma(u_-(s)))^2\right)dxds\\
&&+\int_0^t\int_ {\mathbb{R}^d} \psi(x)k_\epsilon'(v(x,s))(f(u(x,s))-f(u_-(s)))dxds\\
&&+\int_0^t\int_ {\mathbb{R}^d}\psi(x) k_\epsilon'(v(x,s))(\sigma(u(x,s))-\sigma(u_-(s)))dxdW_s.
   \eess
Taking expectation over the above equality and using Lemma \ref{l2.1}, we get
    \bess
\mathbb{E}\Phi_\varepsilon(v(t))
&=&\mathbb{E}\Phi_\varepsilon(v_0)+\mathbb{E}\int_0^t\int_{\mathbb{R}^d}\psi(x) k_\epsilon''(v(x,s))\\
&&\times\left(\frac{1}{2}(\sigma(u(x,s))-\sigma(u_-(s)))^2\right)dxds\\
&&+\mathbb{E}\int_0^t\int_{\mathbb{R}^d}\psi(x) k_\epsilon'(v(x,s))(f(u(x,s))-f(u_-(s)))dxds\\
&\leq&\mathbb{E}\Phi_\varepsilon(v_0)+\frac{L_\sigma}{2}\mathbb{E}\int_0^t\int_{\mathbb{R}^d}\psi(x) k_\epsilon''(v(x,s))
v(x,s)^{2}dxds\\
&&+L_f\mathbb{E}\int_0^t\int_{\mathbb{R}^d}\psi(x)|k_\epsilon'(v(x,s))|v(x,s)|dxds.
  \eess
Note that $\lim\limits_{\varepsilon\rightarrow0}\mathbb{E}\Phi_\varepsilon(v(t))
=\mathbb{E}(\eta(v(t))^2,\psi)$, by
taking the limits termwise as $\varepsilon\rightarrow0$ and using Lemma \ref{l2.1}, we have
    \bess
\mathbb{E}(\eta(v(t))^2,\psi)\leq(L^2_\sigma+2L_f)
\int_0^t\mathbb{E}(\eta(v(s))^2,\psi) ds,
   \eess
which, by means of Gronwall's inequality, implies that
  \bess
\mathbb{E}(\eta(v(t))^2,\psi)=0, \ \ \ \ \forall \ t\in[0,T].
   \eess
Note that for any $R>0$, the above inequality always holds.
It follows from $\psi\geq0$ that $\eta(v(t))=v^-(x,t)=0$ a.s. for a.e. $x\in \mathbb{R}^d$.
which implies that $v^-=0$ a.s. for a.e. $x\in \mathbb{R}^d$, and for any $t\in[0,T]$.

Similarly, if we let  $v=u_+-u(x,t)$, then we can prove that $v^-=0$ a.s. for a.e. $x\in \mathbb{R}^d$, and for any $t\in[0,T]$.
That is to say, (\ref{2.4}) holds. $\Box$

\begin{remark}\lbl{r2.1}
It is remarked that in Theorem \ref{t2.1} the norm of solution of problem
(\ref{2.2}) is maximum in $\mathbb{R}^d$. The advantage of this norm is that
the estimates we obtained hold point-wise. Under this norm, Theorem \ref{t2.1} also holds
if the equation (\ref{2.2}) is
replaced by the equation (\ref{2.6}).
\end{remark}

Now, we consider a special case on a bounded domain $D\subset\mathbb{R}^d$ ($d\geq1$)
  \bes\left\{\begin{array}{lll}
   du=(\Delta u+f(u))dt+udW_t, \qquad t>0,&x\in  D,\\[1.5mm]
   u(x,0)=u_0(x), \ \ \ &x\in D,\\[1.5mm]
   u(x,t)=0, \qquad \qquad \qquad \qquad \qquad t>0,  &x\in\partial D.
    \end{array}\right.\lbl{2.6}\ees
We will establish the relationship
between (\ref{2.6}) and the following SDE
  \bes\left\{\begin{array}{lll}
   dX_t=f(X_t)dt+X_tdW_t, \ \ \qquad t>0,\\[1.5mm]
   X(0)=u_0.
    \end{array}\right.\lbl{2.7}\ees
In order to do that, we
consider the eigenvalue problem for the elliptic equation
  \bess\left\{\begin{array}{llll}
-\Delta \phi=\lambda \phi, \ \ \ \ \ \ \  \ {\rm in} \  D,\\
\phi=0, \ \ \qquad\ \qquad  {\rm on}\ \partial D.
   \end{array}\right. \eess
Then, all the eigenvalues are strictly positive, increasing and
the eigenfunction $\phi_1$ corresponding to the smallest eigenvalue
$\lambda_1$ does not change sign in domain $ D$, as shown in \cite{GTbook}.
Therefore, we   normalize it in such a way that
   \bess
\phi_1(x)\geq0,\ \ \ \ \int_ D \phi_1(x)dx=1.
   \eess
  \begin{theo}\lbl{t2.2} Assume that
     \bess
  (f(u),\phi_1)\leq f((u,\phi_1)).
     \eess
If $0$ is a stable (exponentially mean square stable or stochastically stable)
trivial solution  of (\ref{2.7}), then
$0$ is also a stable (exponentially mean square sable or stochastically stable) trivial solution of (\ref{2.6}),
where we take the spatial norm as $\|u\|_{\phi_1}=(u,\phi_1)$.
   \end{theo}

{\bf Proof.} It follows from \cite[Theorem 2.1]{C2011} and \cite{LD2015} that the solutions
of (\ref{2.6}) keep non-negative almost surely, i.e.,
$u(x,t)\geq 0$, a.s. for almost every $x\in D$ and for all $ t\in[0,T]$.
Moreover, the solutions exist globally.
 Let
    \bess
v(t)=(u,\phi_1)=\int_Du(x,t)\phi_1(x)dx.
   \eess
Then $v$ satisfies
  \bes\left\{\begin{array}{lll}
   dv(t)\leq[-\lambda_1v(t)+f(v(t))]dt+v(t)dW_t, \ \ \qquad t>0,&x\in  D,\\[1.5mm]
   v(0)=v_0=(u_0,\phi_1).
    \end{array}\right.\lbl{2.8}\ees
It is easy to prove that $v(t)$ is a sub-solution of the
following problem
  \bes\left\{\begin{array}{lll}
   dY_t=[-\lambda_1Y_t+f(Y_t)]dt+Y_tdW_t, \ \ \qquad t>0,&x\in  D,\\[1.5mm]
   Y(0)=v_0.
    \end{array}\right.\lbl{2.9}\ees
Since the solutions of (\ref{2.9}) keep non-negative, we obtain that
$v(t)$ is also a sub-solution of (\ref{2.8}).
Set $Z_t=Y_t-v(t)$, similar to the proof of Theorem \ref{t2.1}, we can prove
$Z_t\geq0$ almost surely. Indeed, one can first prove that $Y_t\geq0$ almost surely
by using the same method to the proof of Theorem \ref{t2.1}. Then it follows from the
definition of $v$ that $v(t)=(u,\phi_1)>0$ almost surely. Lastly,
noting that
   \bess
-(Y_t^r-v(t)^r)=-r\xi^{r-1}Z_t\geq0, \ \ {\rm when}\ Z_t\leq0,
   \eess
one can use the same method to the proof of Theorem \ref{t2.1} to get
$Z_t\geq0$ almost surely. Similarly, one can prove $Y_t\leq X_t$ almost
surely.  Therefore, if $0$ is a stable (exponentially mean square stable or stochastically stable)
trivial solution  of (\ref{2.7}) , then
$0$ is also a stable  trivial solution of (\ref{2.6}).
The proof is complete. $\Box$


\medskip

{\bf Example} Let $f(u)=au-ku^r$, where $a\in\mathbb{R},k\geq0,r\geq1$.
$-u^r \geq0$ for $u\leq0$. It follows from \cite[Theorem 2.1]{C2011} and \cite{LD2015} that the solutions
of (\ref{2.6}) with $f(u)=au-ku^r$ keep non-negative almost surely.
By using the H\"{o}lder inequality, we have
   \bess
(u,\phi_1)^r=\left(\int_Du(x,t)\phi_1(x)dx\right)^r\leq\int_Du^r(x,t)\phi_1(x)dx=(u^r,\phi_1).
   \eess
Therefore, all the assumptions of Theorem \ref{t2.2} are satisfied.

\begin{remark}\lbl{r2.2} In particular, if $f(u)=u$, then the stability of the trivial solution
$0$ for (\ref{2.6}) and the following problem
  \bess\left\{\begin{array}{lll}
   dY_t=(1-\lambda_1)Y_tdt+Y_tdW_t, \ \ \qquad t>0,&x\in  D,\\[1.5mm]
   Y(0)=(u_0,\phi_1),
    \end{array}\right. \eess
are equivalent. The above problem is different from (\ref{2.7}) because of the Laplacian operator.
\end{remark}

\section{Bounded domain}
\setcounter{equation}{0}

In this section, we consider the stability results on a bounded domain $D\subset\mathbb{R}^d$.
We focus on the conditions which
induce the solution stable. Meanwhile, we are interested in the difference between the
stochastic partial differential equations and partial differential equations. We will consider
the impact of different noise. Throughout this section, $\|\cdot\|$ means
the norm of $L^2(D)$.

We first consider the following initial boundary problem
 \bes\left\{\begin{array}{lll}
   du=\mu\Delta udt+\sigma dW_t,\ \qquad t>0,&x\in  D,\\[1.5mm]
   u(x,0)=u_0(x), \ \ \ &x\in D,\\[1.5mm]
   u(x,t)=0, \quad\qquad \qquad \qquad  t>0,  &x\in\partial D,
    \end{array}\right.\lbl{0.1}\ees
where $\mu$ and $\sigma$ are positive constants. Then we obtain

\begin{theo}\lbl{t0.1} If $\sigma^2|D|<2\mu\lambda_1\mathbb{E}\|u_0\|^2$, then the
trivial solution $0$ is mean square stable.
\end{theo}

{\bf Proof.} We take the Lyapunov function as $\|u\|^2$. By using It\^{o}'s formula, we have
    \bess
\frac{d}{dt}\mathbb{E}\|u(\cdot,t)\|^2&=&2\mu\mathbb{E}\int_ D u\Delta u(x,t)dx+\sigma^2|D|\nm\\
&=&-2\mu\mathbb{E}\int_ D |\nabla u(x,t)|^2dx+\sigma^2|D|\\
&\leq&-2\mu\lambda_1\mathbb{E}\|u(\cdot,t)\|^2+\sigma^2|D|,
    \eess
where we used the Poincare inequality.
Solving the above inequality, we have
  \bess
\mathbb{E}\|u(\cdot,t)\|^2<\left[\mathbb{E}\|u_0\|^2-\frac{\sigma^2|D|}{2\lambda_1\mu}\right]e^{-2\lambda_1\mu t}
+\frac{\sigma^2|D|}{2\lambda_1\mu}.
   \eess
Note that
    \bess
\left[\mathbb{E}\|u_0\|^2-\frac{\sigma^2|D|}{2\lambda_1\mu}\right]e^{-2\lambda_1\mu t}
+\frac{\sigma^2|D|}{2\lambda_1\mu}\leq \mathbb{E}\|u_0\|^2
   \eess
is equivalent to
   \bess
\left[\mathbb{E}\|u_0\|^2-\frac{\sigma^2|D|}{2\lambda_1\mu}\right]\left[e^{-2\lambda_1\mu t}-1\right]
\leq 0.
   \eess
Therefore, if $\mathbb{E}\|u_0\|^2-\frac{\sigma^2|D|}{2\lambda_1\mu}>0$, that is,
$\sigma^2|D|<2\mu\lambda_1\mathbb{E}\|u_0\|^2$, we have
   \bess
\mathbb{E}\|u(\cdot,t)\|^2\leq\mathbb{E}\|u_0\|^2.
    \eess
The proof is complete. $\Box$

\begin{remark}\lbl{r3.1} It is easy to see that if $\sigma=0$, then the
trivial solution $0$ is stable; and if $\mu=0$, then the trivial solution $0$
will be unstable; and when $\sigma^2>0$, $0$ will be a stable trivial solution under
some more assumptions. Possibly one can say that the additive noise will have a "bad" effect on the stability
of trivial solution.

By using the Chebyshev inequality, one can easily prove the trivia solution $0$ is stochastic
stable without any more assumption, see the next theorem for the proof.
  \end{remark}

We study the impact of additive noise.
Consider the following equation
 \bes\left\{\begin{array}{lll}
   du=(\Delta_p u+f(u))dt+\sigma dW_t, \ \ \qquad t>0,&x\in  D,\\[1.5mm]
   u(x,0)=u_0(x), \ \ \ &x\in D,\\[1.5mm]
   u(x,t)=0, \quad \qquad \qquad \qquad \qquad \qquad t>0,  &x\in\partial D,
    \end{array}\right.\lbl{3.1}\ees
where $\Delta_pu=\nabla\cdot(|\nabla u|^{p-2}\nabla u)$.
Because we only consider the
stability of solutions, throughout this paper we will assume the problem we consider
admits a unique global solution. Let $C_\infty$ be the Sobolev embedding constant satisfying
   \bes
\|u\|_{L^\infty(D)}\leq C_\infty\|u\|_{W^{1,p}(D)},\ \ \ p>d.
    \lbl{3.2}\ees
It is noted that the solution $u$ of (\ref{3.1}) satisfies that
$\|u\|_{W^{1,p}(D)}=\|\nabla u\|_{L^{p}(D)}$.
  \begin{theo}\lbl{t3.1} Assume the nonlinear term satisfies
     \bes
uf(u)\leq au^2+bu^{2m},\ \ \ m\geq1,
   \lbl{3.3}\ees
where $a,b\in\mathbb{R}$. Assume further that $p>\max\{2m,d\}$.

If
   \bes
 a+\frac{2-\gamma}{2}
 \left(\frac{\gamma|D||b|C_\infty^p}{2}\right)^{\frac{\gamma}{2-\gamma}}+\frac{\sigma^2|D|}{2\mathbb{E}\|u_0\|^2}<0,
      \lbl{3.4}\ees
where $\gamma\in(0,2)$ satisfies
    \bess
(2m-2+\gamma)\cdot\frac{2}{\gamma}=p.
   \eess
Then the trivial solution $0$ of (\ref{3.1}) is mean square stable.
  \end{theo}

{\bf Proof.} We remark that when $m=1$, Theorem \ref{t3.1} will become easier, see
\cite{LM1998} for similar results.
We pick a Lyapunov function $V(u)=\|u\|^2$.
By It\^{o}'s formula, taking expectation and integrating with respect to $t$, we have
   \bes
\frac{d}{dt}\mathbb{E}\|u(\cdot,t)\|^2&=&2\mathbb{E}\int_ D u(\Delta_p u(x,t)+f(u))dx+\sigma^2|D|\nm\\
&\leq&-2\mathbb{E}\int_ D |\nabla u(x,t)|^pdx+2\mathbb{E}\int_D(au^2+bu^{2m})dx+\sigma^2|D|.
   \lbl{3.5}\ees
By the Sobolev embedding inequality (\ref{3.2}), we have
    \bes
\|u\|_{L^{2m}}^{2m}=\int_ D |u|^{2m}(x,t)dx&\leq&\|u\|_{L^\infty}^{2m-2}\int_ D u^2(x,t)dx\nm\\
&=&\|u\|_{L^\infty}^{2m-2}\|u\|_{L^2}^{2}\nm\\
&\leq&|D|^{\frac{\gamma}{2}}\|u\|_{L^\infty}^{2m-2+\gamma}\|u\|_{L^2}^{2-\gamma}\nm\\
&\leq&C_\infty^{2m-2+\gamma}|D|^{\frac{\gamma}{2}}\|u\|_{W^{1,p}}^{2m-2+\gamma}\|u\|_{L^2}^{2-\gamma}\nm\\
&\leq&\frac{2-\gamma}{2}
 (C_\infty)^{\frac{2(2m-2+\gamma)}{2-\gamma}}
 \left(\frac{\gamma|b||D|}{2}\right)^{\frac{\gamma}{2-\gamma}}\|u\|_{L^2}^{2}\nm\\
 &&
 +\frac{1}{|b|}\|u\|_{W^{1,p}}^{(2m-2+\gamma)\cdot\frac{2}{\gamma}}.
  \lbl{3.6}\ees
Noting that $p>\max\{2m,d\}$,
there exists a constant $\gamma\in(0,2)$ such that
   \bess
(2m-2+\gamma)\cdot\frac{2}{\gamma}=p.
   \eess
Submitting (\ref{3.6}) into (\ref{3.5}), we have
  \bes
\frac{d}{dt}\mathbb{E}\|u(\cdot,t)\|^2\leq 2\left(a+\hat C\right)\mathbb{E}\|u\|_{L^2}^{2}+\sigma^2|D|,
   \lbl{3.7}\ees
where
   \bess
\hat C=\frac{2-\gamma}{2}
 \left(\frac{\gamma|D||b|C_\infty^p}{2}\right)^{\frac{\gamma}{2-\gamma}}.
    \eess
Solving the differential inequality (\ref{3.7}) gives
   \bes
\mathbb{E}\|u(\cdot,t)\|^2\leq \left(\mathbb{E}\|u_0\|^2+\frac{\sigma^2|D|}{2(a+\hat C)}\right)
e^{2(a+\hat C)t}-\frac{\sigma^2|D|}{2(a+\hat C)}.
  \lbl{3.8}\ees
Note that  $a+\hat C<0$ implies that $e^{2(a+\hat C)t}<1$ for all $t>0$. Furthermore,
the assumption
   \bess
\mathbb{E}\|u_0\|^2+\frac{\sigma^2|D|}{2(a+\hat C)}>0
   \eess
yields that
   \bess
 \mathbb{E}\|u(\cdot,t)\|^2\leq  \mathbb{E} \|u_0\|^2,
    \eess
which completes the proof.
$\Box$

\medskip

{\bf Example} Consider
 \bess\left\{\begin{array}{lll}
   du=(\Delta_4 u-u+\frac{1}{C_\infty^4}u^2)dt+\sigma dW_t, \ \quad t>0,&x\in  (0,1),\\[1.5mm]
   u(x,0)=u_0(x), \ \ \ &x\in (0,1),\\[1.5mm]
   u(x,t)=0, \qquad \qquad \qquad \qquad \qquad \qquad t>0,  &x\in\partial (0,1).
    \end{array}\right.\eess
It is easy to check that $\gamma=1$ satisfies $(2m-2+\gamma)\cdot\frac{2}{\gamma}=p$, where
$m=3/2,p=4$. Then if the initial data satisfies
  \bess
 \frac{\sigma^2}{2\mathbb{E}\|u_0\|^2}<\frac{3}{4},
       \eess
then Theorem \ref{t3.1} shows that the the trivial solution $0$ of the above problem will be mean square stable.

\begin{remark}\lbl{r3.2} (1) In Theorem \ref{t3.1}, we assume the constant $b$ satisfies (\ref{3.4}).
Note that the constant $C_\infty$ depends on the domain $D$, the dimension $d$ and the constant
$p$, and thus it is hard to give a concrete constant in an example. The reason is that
we used the embedding inequality (\ref{3.2}). On the other hand, we can use the following embedding
inequality replaced (\ref{3.2}):
   \bess
W^{1,p}(D)\hookrightarrow W^{1,2m}(D)\hookrightarrow L^{2m}(D), \ \ p>2m.
   \eess
Let $C_{2m,2m}$ be the Sobolev embedding constant, i.e., $\|u\|_{L^{2m}(D)}\leq C_{2m,2m}\|u\|_{W^{1,2m}(D)}$.
Then under the assumptions that $b\leq \frac{2}{C_{2m,2m}}$ and
   \bess
2a\mathbb{E}\|u_0\|^2+\sigma^2|D|<0,
   \eess
the trivial solution $0$ of (\ref{3.1}) is mean square stable.

(2) We now explain why we did not get the results of stochastic stability.
Like the case of stochastic differential equations, we try to use a Lyapunov function $V(u)=\|u\|^{2r}$ with $0<r<1$ to prove
the stochastic stability.
For additive noise, we can not prove that $\|u(\cdot,t)\|^2>0$ for all $t>0$.
In order to use the It\^{o} formula, we consider the following Lyapunov functional
$V(u)=(\|u\|+\kappa)^{2r}$ with $0<r<1$ and $0<\kappa\ll1$.
This leads to the expression
   \bes
\frac{d}{dt}\mathbb{E}(\|u\|^2+\kappa)^{r}&=&2r\mathbb{E}\left[(\|u\|^2+\kappa)^{r-1}\int_ D u(\Delta_p u(x,t)+f(u))dx\right]\nm\\
&&+r\sigma^2|D|\mathbb{E}\left[(\|u\|^2+\kappa)^{r-1}\right]\nm\\
&&
+2\sigma^2r(r-1)\mathbb{E}(\|u\|^2+\kappa)^{r-2}\left(\int_Dudx\right)^2\nm\\
&\leq&\mathbb{E}\left[r(\|u\|^2+\kappa)^{r-1}(a+\hat C)\|u(\cdot,t)\|^{2}\right.\nm\\
&&\left.+r\sigma^2(\|u\|^2+\kappa)^{r-1}\left(|D|+2(r-1)\frac{\left(\int_Dudx\right)^2}{\|u\|^2+\kappa}\right)\right].
   \lbl{3.9}\ees
Due to the difference $\left(\int_Dudx\right)^2$ and $\left(\int_D|u|dx\right)^2$, we can not
get any help to control the term $|D|$. Note that
   \bess
\int_Dudx=0
   \eess
maybe happen, so we can not use this term. Even though the term
$\left(\int_Dudx\right)^2$ is replaced by $\left(\int_D|u|dx\right)^2$, we
can not get the desired result. The reason is the followings.
The H\"{o}lder inequality implies that
   \bess
\int_D|u|dx\leq|D|^{\frac{1}{2}}\left(\int_D|u|^2dx\right)^{\frac{1}{2}}.
   \eess
Consequently,
   \bess
\frac{\|u\|_{L^1}^2}{\|u\|_{L^2}^2}\leq |D|.
  \eess
Hence we can not use the above inequality in (\ref{3.9}). The aim of the
above discussion is to show the last two terms of right-side hand of
(\ref{3.9}) are in the same level, which are different from the first term for
SPDEs.

But for SDEs, there will be another case. In this case, let $|D|=1$, then we have
   \bess
 |D|+2(r-1)\frac{\left(\int_Dudx\right)^2}{\|u\|^2+\kappa}
 =1+\frac{2(r-1)}{1+\kappa}.
    \eess
Taking $0<r<1/2$ such that $2r+\kappa<1$, we get
    \bess
 \frac{d}{dt}\mathbb{E}(\|u\|^2+\kappa)^{r}\leq0.
    \eess
Letting $\kappa\to0$, and using the Chebyshev inequality, we obtain the stochastic stability.
In all, we find there is a significant difference between SPDEs and SDEs in the stability theory.
  \end{remark}

We remark that the noise can be easily generalized the cylindrical Wiener process.
We first generalize the classical results of deterministic reaction-diffusion equation \cite[Theorem 4.2.1, p 166]{YLbook}
to the following equation
 \bes\left\{\begin{array}{lll}
du(x,t)=(\Delta u(x,t)+f(x,t,u))dt+\sigma udW_t, \ \qquad t>0,&x\in  D,\\[1.5mm]
   u(x,0)=u_0(x), \ \ \ &x\in D,\\[1.5mm]
   u(x,t)=0, \qquad\qquad\qquad \qquad \qquad \qquad \qquad \ \ \qquad t>0,  &x\in\partial D.
    \end{array}\right.\lbl{3.10}\ees

  \begin{theo} \lbl{t3.2}
  Assume that $f(x,t,0)=0$, $f\in C^1(D\times[0,\infty)\times(-\infty,\infty))$.

(i) If there exists a constant $\alpha>0$ such that for all $(x,t)\in D\times[0,\infty)$, $\eta\in\mathbb{R}$,
we have
   \bess
f(x,t,\eta)\leq(\lambda_1-\alpha)\eta,
   \eess
then for the initial data $u_0$ satisfying $0\leq u_0(x)\leq\rho \phi_1(x)$ with $\rho>0$, problem (\ref{3.10}) admits
a unique positive solution $u(x,t)$ and the following estimate holds almost surely
   \bess
0\leq u(x,t)\leq\rho e^{-(\alpha+\frac{\sigma^2}{2})t+\sigma W_t}\phi_1(x),\ \ (x,t)\in D\times[0,\infty).
   \eess
Consequently,
   \bes
\mathbb{E}u(x,t) \leq \rho e^{-\alpha t}\phi_1(x).
    \lbl{3.11}\ees
Assume further that the initial data $u_0$ is a deterministic function, we have
   \bes
\mathbb{P}\left\{\int_Du(x,t)dx\leq\rho e^{-(\alpha+\frac{\sigma^2}{2})t}\right\}\geq\frac{1}{2}.
   \lbl{3.12}\ees

(ii) If there exists a constant $\alpha>0$ such that for all $(x,t)\in D\times[0,\infty)$, $\eta\geq0$,
we have
   \bess
f(x,t,\eta)\geq(\lambda_1+\alpha)\eta,
   \eess
then for every $\delta>0$, when $u_0(x)\geq\delta \phi_1(x)$, problem (\ref{3.10}) admits
a unique positive solution $u(x,t)$, which exists globally or finite time blowup. On the lifespan, the following estimate holds almost surely
   \bess
u(x,t)\geq\delta e^{(\alpha -\frac{\sigma^2}{2})t+\sigma W_t}\phi_1(x),\ \ (x,t)\in D\times[0,\infty).
   \eess
Consequently, $\mathbb{E}\|u(t)\|^2\geq \delta e^{\alpha t}\mathbb{E}\|u_0\|^2$.
   \end{theo}

{\bf Proof.} We first change the stochastic reaction-diffusion equation into random
reaction-diffusion, then by using comparison principle, the desired results are obtained.
More precisely, let
$v(x,t)=e^{-\sigma W_t}u(x,t)$, then $v(x,t)$ satisfies that
    \bes\left\{\begin{array}{llll}
\frac{\partial}{\partial t}v(x,t)=\Delta v(x,t)-\frac{\sigma^2}{2}v(x,t)
+ e^{-\sigma W_t}f(x,t,e^{\sigma W_t}v(x,t)), \ \quad & t>0,\ x\in  D,\\[1.5mm]
   v(x,0)=u_0(x), \ \ \ &\qquad \quad x\in D,\\[1.5mm]
   v(x,t)=0, & t>0,  \ x\in\partial D.
    \end{array}\right.\lbl{3.13}\ees
By using the assumptions, we get
   \bess
\frac{\partial}{\partial t}v(x,t)\leq\Delta v(x,t)-\frac{\sigma^2}{2}v(x,t)
+ (\lambda_1-\alpha)v(x,t).
   \eess
It is easy to check that $\bar v=\rho e^{-(\alpha+\frac{\sigma^2}{2}) t}\phi_1(x)$ is an upper solution of (\ref{3.13}) and
$0$ is a lower solution to (\ref{3.13}),
which implies that for $(x,t)\in D\times[0,\infty)$
   \bess
0\leq v(x,t)\leq \rho e^{-(\alpha+\frac{\sigma^2}{2}) t}\phi_1(x)\Longleftrightarrow
0\leq u(x,t)\leq\rho e^{-(\alpha+\frac{\sigma^2}{2})t+\sigma W_t}\phi_1(x), \ \ a.s.,
   \eess
which implies that
  \bess
0\leq\int_Du(x,t)dx\leq\rho e^{-(\alpha+\frac{\sigma^2}{2})t+\sigma W_t}, \ \  a.s..
   \eess
By using $\mathbb{E}[e^{\sigma W_t}]=e^{\frac{\sigma^2}{2}t}$, we have
  \bess
\mathbb{E}u(x,t) \leq e^{\left(\alpha-\lambda_1\right)t}\phi_1(x).
    \eess
Note that
  \bess
\left\{\int_Du(x,t)dx\leq\rho e^{-(\alpha+\frac{\sigma^2}{2})t}\right\}
\Longleftrightarrow\left\{e^{\sigma W_t}\leq1\right\}
\supset\{W_t\leq0\},
   \eess
thus we have
   \bess
\mathbb{P}\left\{\int_Du(x,t)dx\leq\rho e^{-(\alpha+\frac{\sigma^2}{2})t}\right\}
\geq\mathbb{P}\{W_t\leq0\} =\frac{1}{2},
   \eess
which proves (\ref{3.12}).

Next, we prove (ii). Note that
    \bess
\frac{\partial}{\partial t}v(x,t)\geq\Delta v(x,t)-\frac{\sigma^2}{2}v(x,t)
+ (\lambda_1+\alpha)v(x,t).
   \eess
It follows from that $\underline{v}=\delta e^{(\alpha-\frac{\sigma^2}{2})t}$ is a lower solution of
(\ref{3.13}), thus we have the desired inequality. The proof is complete. $\Box$
\begin{remark}\lbl{r3.3} Following Theorem \ref{t3.2}, it is easy to see that
in mean square sense, the solution of (\ref{3.10}) keeps the same properties as
the deterministic case, which is different from the additive noise, see Theorems \ref{t0.1} and \ref{t3.1}.
Of course, the big difference between the stochastic and deterministic
cases is that there exists an event such that whose probability is large than $0$, where the event
is that the solution of stochastic case maybe have exponentially decay. In other words,
in (\ref{3.12}), if $-\frac{\sigma^2}{2}<\alpha<0$, then the solution $u$ of (\ref{3.10}) satisfies
$\|u(t)\|^2\leq \|u_0\|^2e^{\left(\alpha-\lambda_1-\frac{\sigma^2}{2}\right)t}$ with probability $\frac{1}{2}$.
Maybe from here we can say the noise can stabilize the solutions.

The method we used in Theorem \ref{t3.2} is comparison principle, which is different from the
Lyapunov functional method. The inequality (\ref{3.11}) holds pointwise, which is
different from the earlier results. What's more, the index $\alpha-\lambda_1$ is different from
that obtained by Lyapunov method, see the next theorem.
In part (ii) of Theorem \ref{t3.2} implies the unstable condition of the
trivial solution $0$, which is new in this field.

Comparing with the stochastic ordinary different equations, the role of the Laplacian operator in
stochastic reaction-diffusion equations gives a help with $\lambda_1$ in the stability of trivial
solutions. Indeed, the reason is the Poincare inequality.
  \end{remark}

The impact of multiplicative noise in Theorem \ref{t3.2} is
not satisfied. We give the next result.
  \begin{theo}\lbl{t3.3} Assume that $f(x,t,0)=0$ and there exists a constant $K>0$ such that for all $(x,t)\in D\times[0,\infty)$, $uf(x,t,u)\leq Ku^2$. If
$K-\lambda_1+\frac{\sigma^2}{2}\leq0$,
then the trivial solution $0$ is mean square stable and if
   \bes
K-\lambda_1-\frac{\sigma^2}{2}<0,
   \lbl{3.14}\ees
then the trivial solution $0$ is stochastically stable.
  \end{theo}

{\bf Proof.}
Taking Lyapunov function $V(u)=\|u\|^2$, we have
   \bess
\frac{d}{dt}\mathbb{E}\|u(t)\|^2&=&-\int_D|\nabla u(x,t)|^2dx+\frac{\sigma^2}{2}\mathbb{E}\|u(t)\|^2+\mathbb{E}\int_Duf(x,t,u)dx\\
&\leq&\left(K-\lambda_1+\frac{\sigma^2}{2}\right)\mathbb{E}\|u(t)\|^2,
   \eess
which yields that
   \bess
\mathbb{E}\|u(t)\|^2\leq \mathbb{E}\|u_0\|^2e^{\left(K-\lambda_1+\frac{\sigma^2}{2}\right)t}.
    \eess

Next, we use a Lyapunov function $V(u)=\|u\|^{2r}$ with $0<r<1$ to prove
the stochastic stability. Note that in Theorem
\ref{t3.2}, we proved that the solution $u\geq0$ almost surely and thus
we can choose $\|u\|^{2r}$  as Lyapunov functional. This leads to the expression
   \bes
\frac{d}{dt}\mathbb{E}\|u(\cdot,t)\|^{2r}&=&2r\mathbb{E}\left[\|u(\cdot,t)\|^{2r-2}\int_ D u(\Delta u(x,t)+f(x,t,u))dx\right]\nm\\
&&+r\sigma^2\mathbb{E}\|u(\cdot,t)\|^{2r}
+2\sigma^2r(r-1)\mathbb{E}\|u(\cdot,t)\|^{2r}\nm\\
&\leq&\mathbb{E}\left[2r\|u(\cdot,t)\|^{2r}\left(K-\lambda_1+\sigma^2r-\frac{\sigma^2}{2}\right)\right].
   \lbl{3.15}\ees
If  $K-\lambda_1-\frac{\sigma^2}{2}<0$, we can choose $0<r<\frac{1}{2}$ such that
    \bess
K-\lambda_1+\sigma^2r-\frac{\sigma^2}{2}\leq0,
   \eess
then from the Chebyshev inequality, stochastic stability for the solution
of (\ref{3.1}) follows from (\ref{3.15}). The proof is complete. $\Box$

\begin{remark}\lbl{r3.4} It follows from Theorem \ref{t3.3} that the multiplicative
noise can make solution stable in sense of stochastically stable. Comparing Theorem \ref{t3.3} with
\ref{t3.2}, we can take $K=\lambda_1+\sigma^2/2-\varepsilon>\lambda_1-\alpha$, where $0<\varepsilon\ll1$.
\end{remark}

In the above Theorems, we assume that the noise term satisfies the global Lipschitz condition.
In the following theorem, we will see that the assumption can be weaken as local Lipschitz condition.
In order to do this, we consider the following equation
      \bes\left\{\begin{array}{lll}
   du=(\Delta u-k_1u^{r})dt+k_2u^mdW_t(x), \ \ \qquad &t>0,\ x\in  D,\\[1.5mm]
   u(x,0)=u_0(x), \ \ \ &\qquad\ \ \ x\in D,\\[1.5mm]
   u(x,t)=0, \ \ \ \ &t>0,\ x\in\partial D,
    \end{array}\right.\lbl{3.17}\ees
where $k_1,k_2,r$ and $m$ are positive constants.  Moreover,
$W_t(x)$ is a Wiener random field with the covariance function $q$.
In our paper \cite{LD2015}, we proved the existence of global solution of (\ref{3.17}) under the assumptions
of Theorem \ref{t3.4}. Moreover, we proved the solutions keep non-negative almost surely.

  \begin{theo}\lbl{t3.4} Assume that $r$ is an odd number and $1<m<\frac{1+r}{2}$.
Assume further that there exists a positive constant $q_0$ such that the covariance function $q(x,y)$ satisfies the condition $\sup_{x,y\in\bar D}q(x,y)\leq q_0$.
Then if $\hat\lambda<\lambda_1$, then the trivial solution $0$ is exponentially mean square stable with the index $\lambda_1-\hat\lambda$, where
   \bess
\hat\lambda:=\frac{r+1-2m}{r-1}\left(\frac{k_1(r-1)}{2m-2}\right)^{-\frac{2m-2}{r+1-2m}}(q_0k_2)^{\frac{r-1}{r+1-2m}}.
  \eess
In particular, when $m=2$ and $r>3$, we
assume further that there exists a positive constant $q_1$ such that
the covariance function $q(x,y)$ satisfies the condition $\sup_{x,y\in\bar D}q(x,y)\geq q_1$.
If
   \bess
\lambda_1>\frac{r-3}{r-1}\left(\frac{k_1(r-1)}{2}\right)^{-\frac{2}{r-3}}(q_0k_2)^{\frac{r-1}{r-3}}-2k_2^2q_1,
   \eess
then the trivial solution $0$ is stochastic stable.
  \end{theo}

{\bf Proof.} The proof is similar to that of Theorem \ref{t3.3} and we give outline of the
proof for completeness.
Taking Lyapunov function $V(u)=\|u\|^2$, we have the following inequality
   \bess
\frac{d}{dt}\mathbb{E}\|u(t)\|^2&=&-2\int_D|\nabla u(x,t)|^2dx+k_2^2\mathbb{E}\int_Du^{2m}(x,t)q(x,x)dx-2k_1\mathbb{E}\int_Du^{r+1}(x,t)dx\\
&\leq&-2\lambda_1\mathbb{E}\|u(t)\|^2
-2k_1\mathbb{E}\|u(t)\|_{L^{1+r}}^{1+r}+k_2^2q_0\mathbb{E}\|u(t)\|_{L^{2m}}^{2m}.
   \eess
By using the interpolation inequality
  \bess
\|u\|_{L^r}\leq\|u\|^\theta_{L^p}\|u\|^{1-\theta}_{L^q},
 \eess
with $r=2m,\,p=2$ and $q=4$, we have
  \bess
q_0\|u\|_{L^{2m}}^{2m}&\leq&q_0\|u\|_{L^2}^{2m\theta}\|u\|_{L^{1+r}}^{2m(1-\theta)}\nm\\
&\leq&k_1\|u\|_{L^{1+r}}^{2m(1-\theta)\frac{1}{1-m\theta}}+
m\theta\left(\frac{k_1}{1-m\theta}\right)^{-\frac{1-m\theta}{m\theta}}
(q_0k_2)^{\frac{1}{m\theta}}\|u\|_{L^2}^{2}\nm\\
&=&k_1\|u\|_{L^{1+r}}^{1+r}+\hat\lambda\|u\|_{L^2}^{2},
    \eess
 where
   \bess
 \theta=\frac{r+1-2m}{mr-m},\ \  \hat\lambda:=\frac{r+1-2m}{r-1}\left(\frac{k_1(r-1)}{2m-2}\right)^{-\frac{2m-2}{r+1-2m}}(q_0k_2)^{\frac{r-1}{r+1-2m}}.
  \eess
Consequently, we have
   \bess
\mathbb{E}\|u(t)\|^2\leq \mathbb{E}\|u_0\|^2e^{-\left(\lambda_1-\hat\lambda\right)t},
    \eess
where implies that the trivial solution $0$ is exponentially mean square stable.

Next, we use a Lyapunov function $V(u)=\|u\|^{2\gamma}$ with $0<\gamma<1$ to prove
the stochastic stability (note that the solutions of (\ref{3.17}) is non-negative function, see \cite{LD2015}).
This leads to the expression
   \bess
\frac{d}{dt}\mathbb{E}\|u(\cdot,t)\|^{2\gamma}&=&\gamma\mathbb{E}
\left[\|u(\cdot,t)\|^{2\gamma-2}\int_ D u(\Delta u(x,t)-k_1u^r(x,t))dx\right]\nm\\
&&+\gamma k_2^2\mathbb{E}\|u(\cdot,t)\|^{2\gamma-2}\int_Du^{2m}dx\nm\\
&&
+2k_2^2\gamma(\gamma-1)\mathbb{E}\|u(\cdot,t)\|^{2\gamma-4}\int_D\int_Dq(x,y)u^{m}(x,t)
u^m(y,t)dx\\
&\leq&\gamma\mathbb{E}
\left[\|u(\cdot,t)\|^{2\gamma-2}\int_ D u(\Delta u(x,t)-k_1u^r(x,t))dx\right]\nm\\
&&+\gamma k_2^2\mathbb{E}\|u(\cdot,t)\|^{2\gamma-2}\int_Du^{2m}dx
+2k_2^2q_0\gamma(\gamma-1)\mathbb{E}\|u(\cdot,t)\|^{2\gamma}.
   \eess
Then by using similar method in proving mean square stable, we can choose
$0<r<1$ such that
   \bess
\frac{d}{dt}\mathbb{E}\|u(\cdot,t)\|^{2\gamma}\leq0.
   \eess
From the Chebyshev inequality, stochastic stability for the solution
of (\ref{3.17}) is obtained. The proof is complete. $\Box$

In  paper \cite{TWW2020}, the authors considered the following problem
      \bes\left\{\begin{array}{lll}
   dX_t=X_t(b(X_t)+k_1X^{m-1}_t)dt+k_2X^{\frac{m+1}{2}}_t\phi(X_t)dW_t, \ \ \qquad &t>0,\ \\[1.5mm]
  X_0=x>0,
    \end{array}\right.\lbl{3.18}\ees
where $k_1,k_2\in\mathbb{R},\,m\geq1$. In \cite{LDWW2018}, we considered
the competition between the nonlinear term and noise term. The
result of  \cite{TWW2020} generalized the results of \cite{LDWW2018}.
Now we first recall the main results of  \cite{TWW2020}.

 \begin{prop}\lbl{p3.1}\cite[Theorem 1.1]{TWW2020} Let $k_1$ be a real number which is not zero.
 Assume $rb(r)\in C^1(\mathbb{R}_+)$ and there exist two positive
 numbers $c_0$ and $m_0(<m)$ such that
    \bess
 |b(r)|\leq c_0(1+r^{m_0-1}),\ \ \ r\in\mathbb{R}_+.
    \eess
 Assume in addition that $r^{\frac{m+1}{2}}\phi(r)\in C^1(\mathbb{R}_+)$. Let $\beta\in(0,1)$ and suppose there
 is a $r_0>0$ such that
    \bess
 \inf_{r\geq r_0}\phi(r)>\sqrt{\frac{2|k_1|}{(1-\beta)k_2^2}}.
   \eess
There is a unique solution $X_t(x)$ for (\ref{3.18}) on $t\geq0$
and the solution is positive for all $t\geq0$ almost surely. Moreover,
for every $T>0$
    \bess
\sup_{0\leq t\leq T}\mathbb{E}X_t^\beta(x)<+\infty.
   \eess
\end{prop}

Moreover, \cite[Theorem 1.2]{TWW2020} shows that the result in
proposition is sharp. More precisely, if there exists $\gamma\in(\beta,1)$ such that
    \bess
 \sup_{r\geq r_0}\phi(r)<\sqrt{\frac{2|k_1|}{(1-\gamma)k_2^2}},
   \eess
then there is a real number $T_0>0$ such that
    \bess
\sup_{0\leq t\leq T}\mathbb{E}X_t^\gamma(x)=+\infty.
   \eess
The above results implies that the
trivial solution $0$ is not mean square stable. But
the stochastic stability would be possible. In the following
result, we will give a positive answer.

\begin{theo}\lbl{t3.6} Let all the assumptions of Proposition \ref{p3.1} hold.
Assume further that
    \bess
rb(r)\leq c_1r+c_2r^{m_0},\ \ 1<m_0<m,\ \ r\in\mathbb{R}_+,\ c_1<0<c_2.
    \eess

(i) If $ \inf_{r\geq 0}\phi(r)\geq1$ and
    \bess
c_1+\frac{[p(k_1-\frac{k_2^2}{2})]^{-p/q}}{q}<0,
   \eess
where
   \bess
p=\frac{m-1}{m_0-1}, \ \ \ q=\frac{m-1}{m-m_0}.
   \eess
Then the trivial solution $0$ is stochastic stable.

(ii) If $\phi(r)=r^{\frac{\alpha}{2}}$ with $\alpha\geq0$ and
    \bess
&&\frac{m-1}{\alpha+m-1}\left(\frac{-c_1(\alpha+m-1)}{2\alpha}\right)^{-\frac{\alpha}{m-1}}k_1^{\frac{\alpha+m-1}{m-1}}\\
&&
+\frac{m_0-1}{\alpha+m-1}\left(\frac{-c_1(\alpha+m-1)}{2(\alpha+m-m_0)}\right)^{-\frac{\alpha+m-m_0}{m_0-1}}
c_2^{\frac{\alpha+m-1}{m_0-1}}<\frac{k_2^2}{2}.
   \eess
Then the trivial solution $0$ is stochastic stable.
 \end{theo}

{\bf Proof.}
(i) We use a Lyapunov function $V(X)=|X|^{\beta}$ with $0<\beta<1$ being
fixed later to prove
the stochastic stability (noting that the solutions is positive almost surely, or
one can use $(|X|+\kappa)^{\beta}$ to replace $|X|^{\beta}$ and then let
$\kappa\to0$). This leads to the expression
   \bess
\frac{d}{dt}\mathbb{E}|X|^{\beta}&=&\beta\mathbb{E}
|X|^{\beta}(b(X_t)+k_1X^{m-1}_t)+\frac{1}{2}k_2^2\beta(\beta-1)\mathbb{E}|X_t|^{\beta+m-1}\phi^2(X_t)\\
&\leq&\beta\mathbb{E}|X|^{\beta}\left[c_1+c_2X_t^{m_0-1}+\left(k_1+(\beta-1)\frac{k_2^2}{2}\right)
|X_t|^{m-1}\right].
    \eess
Set $0<\beta\ll1$ such that
   \bess
c_1+\frac{[p(k_1-\frac{(1-\beta)k_2^2}{2})]^{-p/q}}{q}\leq0.
   \eess
By using the $\varepsilon$-Young inequality, we have
   \bess
\frac{d}{dt}\mathbb{E}|X|^{\beta}\leq0.
    \eess
From the Chebyshev inequality, stochastic stability for the solution
of (\ref{3.18}) is obtained.

(ii) In this case: $\phi(r)=r^{\frac{\alpha}{2}}$ with $\alpha\geq0$. For every $\beta\in(0,1)$,
if we take
$r_0=\left(\sqrt{\frac{3|k_1|}{k_2(1-\beta)}}\right)^{\frac{2}{\alpha}}$, then
   \bess
\inf_{r\geq r_0}\phi(r)=\inf_{r\geq r_0}|r|^{\alpha/2}=\sqrt{\frac{3|k_1|}{k_2(1-\beta)}}
>\sqrt{\frac{2|k_1|}{(1-\beta)k_2^2}}.
   \eess
Hence Proposition \ref{p3.1} holds for this case.

Similar to case (i), we use a Lyapunov function $V(X)=|X|^{\beta}$ with $0<\beta<1$.
This leads to the expression
   \bess
\frac{d}{dt}\mathbb{E}|X|^{\beta}&=&\beta\mathbb{E}
|X|^{\beta}(b(X_t)+k_1X^{m-1}_t)+\frac{1}{2}k_2^2\beta(\beta-1)\mathbb{E}|X_t|^{\beta+\alpha+m-1}\\
&\leq&\beta\mathbb{E}|X|^{\beta}\left[c_1+c_2X_t^{m_0-1}+k_1|X_t|^{m-1}+(\beta-1)\frac{k_2^2}{2}
|X_t|^{\alpha+m-1}\right].
    \eess
Set $0<\beta\ll1$ such that
\bess
&&\frac{m-1}{\alpha+m-1}\left(\frac{-c_1(\alpha+m-1)}{2\alpha}\right)^{-\frac{\alpha}{m-1}}k_1^{\frac{\alpha+m-1}{m-1}}\\
&&
+\frac{m_0-1}{\alpha+m-1}\left(\frac{-c_1(\alpha+m-1)}{2(\alpha+m-m_0)}\right)^{-\frac{\alpha+m-m_0}{m_0-1}}
c_2^{\frac{\alpha+m-1}{m_0-1}}\leq\frac{k_2^2}{2}(1-\beta).
   \eess
By using the $\varepsilon$-Young inequality, we have
   \bess
\frac{d}{dt}\mathbb{E}|X|^{\beta}\leq0.
    \eess
From the Chebyshev inequality, stochastic stability for the solution
of (\ref{3.18}) is obtained. $\Box$

\begin{remark}\lbl{r3.5}
Theorem \ref{t3.6} is new for SDEs. When $k_2=0$, the solution of
(\ref{3.18}) will blow up in finite time, and thus  Theorem \ref{t3.6}
implies that the multiplicative noise can make the solution stable.
\end{remark}

Unfortunately, for SPDEs, we can not get the similar result to
Theorem \ref{t3.6}. Before we end this section, we give the reason.
For simplicity, we consider the following problem
      \bes\left\{\begin{array}{lll}
   du=(\Delta u+k_1u^{r})dt+k_2u^mdW_t, \ \ \qquad &t>0,\ x\in  D,\\[1.5mm]
   u(x,0)=u_0(x), \ \ \ &\qquad\ \ \ x\in D,\\[1.5mm]
   u(x,t)=0, \ \ \ \ &t>0,\ x\in\partial D,
    \end{array}\right.\lbl{3.19}\ees
where $k_1,k_2,r,m\in\mathbb{R}$. Under the condition $r<m$, the
existence of global solution was established by \cite{LW2020}. In the
following, we set forth the reason why we can not get that the trivial solution
$0$ is stochastic stable.

We use a Lyapunov function $V(u)=\|u\|^{2\gamma}$ with $0<\gamma<1/2$ to prove
the stochastic stability (if we worry about $\|u\|=0$, we can use
$(\|u\|^2+\kappa)^\gamma$ instead and then let $\kappa\to0$). This leads to the expression
    \bes
\frac{d}{dt}\mathbb{E}\|u(\cdot,t)\|^{2\gamma}&=&\gamma\mathbb{E}
\left[\|u(\cdot,t)\|^{2\gamma-2}\int_ D u(\Delta u(x,t)+k_1u^r(x,t))dx\right]\nm\\
&&+\gamma\mathbb{E}\left[\|u(\cdot,t)\|^{2\gamma-2}\int_ Dq(x,x)u^{2m}dx\right]\nm\\
&&
+2k_2^2\gamma(\gamma-1)\mathbb{E}\|u(\cdot,t)\|^{2\gamma-4}\int_D\int_Dq(x,y)u^{m+1}(x,t)
u^{m+1}(y,t)dxdy\nm\\
&\leq&-\gamma\lambda_1\mathbb{E}\|u(\cdot,t)\|^{2\gamma}+
k_1\gamma\mathbb{E}\left[\|u(\cdot,t)\|^{2\gamma-2}
\|u(\cdot,t)\|_{L^{1+r}}^{1+r}\right]\nm\\
&&+\gamma\mathbb{E}\left[\|u(\cdot,t)\|^{2\gamma-2}\left(q_0\|u\|_{L^{2m}}^{2m}
+
2q_1k_2^2(\gamma-1)\frac{\|u(\cdot,t)\|_{L^{m+1}}^{2(m+1)}}{\|u(\cdot,t)\|^{2}}\right)\right].
    \lbl{3.20}\ees
H\"{o}lder inequality implies that
   \bess
\int_D| u|^{m+1}(x,t)dx \leq \left(\int_D| u|^{2m}(x,t)dx\right)^{\frac{1}{2}}\left(\int_D| u|^{2}(x,t)dx\right)^{\frac{1}{2}}.
   \eess
For the last term of right-side hand of (\ref{3.20}), we want to get
  \bess
\frac{\|u(\cdot,t)\|_{L^{m}}^{2m}}{\|u(\cdot,t)\|^{2}}\geq\|u\|_{L^{2m}}^{2m},
   \eess
which is a contradiction with respect to the above H\"{o}lder inequality.
This shows no difference from Remark \ref{r3.2}.

\section{Whole space}
\setcounter{equation}{0}

In this section, we consider the impact of noise in the whole space.
 \bes\left\{\begin{array}{lll}
   du=(\Delta u+f(t,u))dt+\sigma(t,u)dW(x,t) \ \ \qquad t>0,&x\in  \mathbb{R},\\[1.5mm]
   u(x,0)=u_0(x), \ \ \ &x\in \mathbb{R},
    \end{array}\right.\lbl{4.1}\ees
where $W(x,t)$ is the space-time white noise. Throughout this section, we assume that $f(t,0)=\sigma(t,0)=0$. A mild solution to (\ref{4.1}) in sense of Walsh \cite{walsh1986} is any
$u$ which is adapted to the filtration generated by the
white noise and satisfies the following evolution equation
   \bess
u(x,t)=\int_{\mathbb{R}}K(x-y,t)u_0(y)dy
+\int_0^t\int_{\mathbb{R}}K(x-y,t-s)\sigma(u,y,s)dW(y,s)dy,
   \eess
where $K(x,t)$ denotes the heat kernel of Laplacian operator.

  \begin{theo}\lbl{t4.1} Assume that there exist non-negative constants $\alpha$, $\beta$ and $\gamma$ such that
  \bes
|f(u)+\alpha u|\leq\beta(t) |u|, \ \ |\sigma(t,u)|\leq\gamma(t) |u|,\ \
2\int_0^t\beta(s)ds +\frac{1}{2\sqrt{\pi}}\int_0^t\frac{\gamma(s)}{(t-s)^{-\frac{1}{2}}}ds<1.
  \lbl{4.1}\ees
Then the trivial solution $0$ is exponentially mean square stable with index $2\alpha$.
  \end{theo}

{\bf Proof.} Note that $f$ and $\sigma$ satisfy the global Lipschitz condition,
one can prove that (\ref{4.1}) admits a unique global solution.
Let $v(x,t)=e^{\alpha t}u(x,t)$. Then $v$ satisfies
 \bes\left\{\begin{array}{lll}
   dv=(\Delta v+f(e^{-\alpha t}v)+\alpha e^{-\alpha t}v)dt+\sigma(t,e^{-\alpha t}v)dW(x,t), \ \qquad t>0,&x\in  \mathbb{R},\\[1.5mm]
   v(x,0)=u_0(x), \ \ \ &x\in \mathbb{R}.
    \end{array}\right.\lbl{4.2}\ees
By taking the second moment and using the Walsh isometry,
we get for $t\in[0,T]$
   \bess
\mathbb{E}|v(x,t)|^2&=&\left(\int_{\mathbb{R}}K(x-y,t)u_0(y)dy+\int_0^t\int_{\mathbb{R}}K(x-y,t-s)
[f(e^{-\alpha s}v)+\alpha e^{-\alpha s}v]dyds\right)^2\\
&&
+\int_0^t\int_{\mathbb{R}}K^2(x-y,t-s)\mathbb{E}\sigma^2(s,e^{-\alpha t}v)dyds\\
&\leq&2\max_{x\in\mathbb{R}}|u_0|^2(x)+\left(2\int_0^t\beta(s)ds +\frac{1}{2\sqrt{\pi}}\int_0^t\frac{\gamma(s)}{(t-s)^{-\frac{1}{2}}}ds\right)\max_{(x,t)\in\mathbb{R}\times[0,T]}
\mathbb{E}v^2(x,t),
   \eess
which implies that
   \bess
\mathbb{E}|v(x,t)|^2\leq C\max_{x\in\mathbb{R}}|u_0|^2(x),
   \eess
where $C$ is independent of $t$. And thus we complete the proof. $\Box$

The reason why we used the properties of heat kernel is that we did not know the existence of local strong solution of problem (\ref{4.1}) with space-time white noise. For general $d$, we have the following result.

\begin{theo} \lbl{t4.2} Consider the Cauchy problem
\bes\left\{\begin{array}{lll}
du=(\Delta u+f(t,u))dt+\sigma(t,u)dW_t \ \ \qquad t>0,&x\in  \mathbb{R}^d,\\[1.5mm]
u(x,0)=u_0(x), \ \ \ &x\in \mathbb{R}^d.
\end{array}\right.\lbl{4.4}\ees
Assume that there exists a constant $K>0$ such that for all $(x,t)\in \mathbb{R}^d\times[0,\infty)$,
    \bess
uf(t,u)\leq Ku^2,\ \ \ \sigma^2(t,u)\leq \sigma u^2, \ \ \sigma>0.
   \eess
Then if
$K+\frac{\sigma^2}{2}\leq0$,
then the trivial solution $0$ is mean square stable and if
$K-\frac{\sigma^2}{2}<0$,
then the trivial solution $0$ is stochastically stable.
  \end{theo}

The proof of Theorem \ref{t4.2} is similar to Theorem \ref{t3.3} and we omit it here.
Meanwhile, we remark that in the whole space the operator $\Delta$ will have no help to
the stability of trivial solution.

\bigskip

\noindent {\bf Acknowledgment} This research was partly supported by the NSF of China grants 11771123, 11501577, 11626085.

 \end{document}